\documentclass{article}
\usepackage{amsmath}
\usepackage{amssymb}
\begin{document}
\def\eq#1{{\rm(\ref{#1})}}
\newtheorem{thm}{Theorem}[section]
\newtheorem{prop}[thm]{Proposition}
\newenvironment{dfn}{\medskip\refstepcounter{thm}
\noindent{\bf Definition \thesection.\arabic{thm}\ }}{\medskip}
\newenvironment{ex}{\medskip\refstepcounter{thm}
\noindent{\bf Example \thesection.\arabic{thm}\ }}{\medskip}
\def\dim{\mathop{\rm dim}}
\def\Re{\mathop{\rm Re}}
\def\Im{\mathop{\rm Im}}
\def\Hol{{\textstyle\mathop{\rm Hol}}}
\def\vol{\mathop{\rm vol}}
\def\U{\mathbin{\rm U}}
\def\SU{\mathop{\rm SU}}
\def\sech{{\textstyle\mathop{\rm sech}}}
\def\ge{\geqslant} 
\def\le{\leqslant} 
\def\R{\mathbin{\mathbb R}}
\def\C{\mathbin{\mathbb C}}
\def\al{\alpha}
\def\be{\beta}
\def\la{\lambda}
\def\ga{\gamma}
\def\de{\delta}
\def\ep{\epsilon}
\def\om{\omega}
\def\th{\theta}
\def\vp{\varphi}
\def\De{\Delta}
\def\La{\Lambda}
\def\Om{\Omega}
\def\d{{\rm d}}
\def\pd{\partial}
\def\bs{\boldsymbol}
\def\w{\wedge}
\def\sm{\setminus}
\def\iy{\infty}
\def\ra{\rightarrow}
\def\t{\times}
\def\ha{{\textstyle\frac{1}{2}}}
\def\ti{\tilde}
\def\ms#1{\vert#1\vert^2}
\def\md#1{\vert #1 \vert}
\def\bmd#1{\big\vert #1 \big\vert}
\title{$\U(1)$-invariant special Lagrangian 3-folds in $\C^3$ \\
and special Lagrangian fibrations}
\author{Dominic Joyce \\ Lincoln College, Oxford}
\date{}
\maketitle

\section{Introduction}
\label{uf1}

Special Lagrangian submanifolds (SL $m$-folds) are a distinguished
class of real $m$-dimensional minimal submanifolds in $\C^m$, which are
calibrated with respect to the $m$-form $\Re(\d z_1\w\cdots\w\d z_m)$.
They can also be defined in (almost) Calabi--Yau manifolds, are
important in String Theory, and are expected to play a r\^ole in the
eventual explanation of Mirror Symmetry between Calabi--Yau 3-folds.

This paper surveys three papers \cite{Joyc3,Joyc4,Joyc5} studying
special Lagrangian 3-folds $N$ in $\C^3$ invariant under the
$\U(1)$-action
\begin{equation}                                             
{\rm e}^{i\th}:(z_1,z_2,z_3)\mapsto
({\rm e}^{i\th}z_1,{\rm e}^{-i\th}z_2,z_3)
\quad\text{for ${\rm e}^{i\th}\in\U(1)$,}
\label{uf1eq1}
\end{equation}
and also the sequel \cite{Joyc6}, which applies their results to study
the {\it SYZ Conjecture} about Mirror Symmetry of Calabi--Yau 3-folds.

Locally we can write $N$ in the form
\begin{equation}
\begin{split}
N=\Bigl\{(z_1&,z_2,z_3)\in\C^3:
\Im(z_3)=u\bigl(\Re(z_3),\Im(z_1z_2)\bigr),\\
&\Re(z_1z_2)=v\bigl(\Re(z_3),\Im(z_1z_2)\bigr),
\quad\ms{z_1}-\ms{z_2}=2a\Bigr\},
\end{split}
\label{uf1eq2}
\end{equation}
where $a\in\R$ and $u,v:\R^2\ra\R$ are continuous functions. Then
$N$ is an SL 3-fold if and only if $u,v$ satisfy
\begin{equation}
\frac{\pd u}{\pd x}=\frac{\pd v}{\pd y}\quad\text{and}\quad
\frac{\pd v}{\pd x}=-2\bigl(v^2+y^2+a^2\bigr)^{1/2}\frac{\pd u}{\pd y}.
\label{uf1eq3}
\end{equation}
If $a\ne 0$ then \eq{uf1eq3} is {\it elliptic}, and so solutions $u,v$
of \eq{uf1eq3} are automatically real analytic, and the corresponding
SL 3-folds $N$ are nonsingular.

However, if $a=0$ then at points $(x,0)$ with $v(x,0)=0$ the factor
$(v^2+y^2+a^2)^{1/2}$ becomes zero, and \eq{uf1eq3} is no longer
elliptic. Because of this, when $a=0$ the appropriate thing to do
is consider {\it weak\/} solutions of \eq{uf1eq3}, which may have
{\it singular points} $(x,0)$ with $v(x,0)=0$. At such a point
$u,v$ may not be differentiable, and $\bigl(0,0,x+iu(x,0)\bigr)$
is a singular point of the SL 3-fold~$N$.

Equation \eq{uf1eq3} is a {\it nonlinear Cauchy--Riemann equation}, so
that if $u,v$ is a solution then $u+iv$ is a bit like a holomorphic
function of $x+iy$. Therefore we may use ideas and methods from
complex analysis to study the solutions of \eq{uf1eq3}, and the
corresponding SL 3-folds~$N$.

Section \ref{uf2} introduces special Lagrangian geometry, and
\S\ref{uf3} recalls some analytic background. Section \ref{uf4}
begins the discussion of SL 3-folds of the form \eq{uf1eq2},
looking at equation \eq{uf1eq3} and what properties we expect
of singular solutions. Some examples of solutions $u,v$ to
\eq{uf1eq3} are given in~\S\ref{uf5}.

Section \ref{uf6} rewrites \eq{uf1eq3} in terms of a {\it potential}
$f$ with $\frac{\pd f}{\pd y}=u$ and $\frac{\pd f}{\pd x}=v$. Then
$f$ satisfies the equation
\begin{equation}
\Bigl(\Bigl(\frac{\pd f}{\pd x}\Bigr)^2+y^2+a^2\Bigr)^{-1/2}
\frac{\pd^2f}{\pd x^2}+2\,\frac{\pd^2f}{\pd y^2}=0.
\label{uf1eq4}
\end{equation}
We can prove existence and uniqueness for the Dirichlet problem
for \eq{uf1eq4} on a suitable class of convex domains in $\R^2$.
This yields existence and uniqueness results for $\U(1)$-invariant
SL 3-folds $N$ in $\C^3$ with boundary conditions, including
singular SL 3-folds.

Section \ref{uf7} studies zeroes of $(u_1,v_1)-(u_2,v_2)$ when $u_j,
v_j$ satisfy \eq{uf1eq3} for $j=1,2$. As $(u_1,v_1)-(u_2,v_2)$ acts
like a holomorphic function, if $(u_1,v_1)\not\equiv(u_2,v_2)$ zeroes
are isolated, and have a positive integer {\it multiplicity}. This
is applied in \S\ref{uf8} to singular solutions $u,v$ of \eq{uf1eq3}
with $a=0$. We find that either $u(x,-y)\equiv u(x,y)$ and $v(x,-y)
\equiv -v(x,y)$, so that $u,v$ is singular all along the $x$-axis,
or else singularities are isolated, with a positive integer
{\it multiplicity} and one of two {\it types}.

Section \ref{uf9} uses the material of \S\ref{uf6}--\S\ref{uf7} to
construct large families of $\U(1)$-invariant {\it special Lagrangian
fibrations} of open subsets of $\C^3$. Finally, \S\ref{uf10} discusses
the {\it SYZ Conjecture}, and summarizes the conclusions of \cite{Joyc6}
on special Lagrangian fibrations of (almost) Calabi--Yau 3-folds.
\medskip

\noindent{\it Acknowledgements.} I would like to thank Robert Bryant
and Mark Gross for helpful conversations. I was supported by an EPSRC
Advanced Fellowship whilst writing this paper. This article was written
for the Proceedings of the $9^{\rm th}$ G\"okova Geometry-Topology
Conference, 2002, supported by TUBITAK, and I am grateful to the
organizers for inviting me.

\section{Special Lagrangian geometry}
\label{uf2}

We shall define special Lagrangian submanifolds first in $\C^m$
and then in {\it almost Calabi--Yau manifolds}, a generalization
of Calabi--Yau manifolds. For introductions to special Lagrangian
geometry, see Harvey and Lawson \cite[\S III]{HaLa} or the
author~\cite{Joyc2}.

\subsection{Special Lagrangian submanifolds in $\C^m$}
\label{uf21}

We begin by defining {\it calibrated submanifolds}, following
Harvey and Lawson~\cite{HaLa}.

\begin{dfn} Let $(M,g)$ be a Riemannian manifold. An {\it oriented
tangent $k$-plane} $V$ on $M$ is a vector subspace $V$ of
some tangent space $T_xM$ to $M$ with $\dim V=k$, equipped
with an orientation. If $V$ is an oriented tangent $k$-plane
on $M$ then $g\vert_V$ is a Euclidean metric on $V$, so 
combining $g\vert_V$ with the orientation on $V$ gives a 
natural {\it volume form} $\vol_V$ on $V$, which is a 
$k$-form on~$V$.

Now let $\vp$ be a closed $k$-form on $M$. We say that
$\vp$ is a {\it calibration} on $M$ if for every oriented
$k$-plane $V$ on $M$ we have $\vp\vert_V\le \vol_V$. Here
$\vp\vert_V=\al\cdot\vol_V$ for some $\al\in\R$, and 
$\vp\vert_V\le\vol_V$ if $\al\le 1$. Let $N$ be an 
oriented submanifold of $M$ with dimension $k$. Then 
each tangent space $T_xN$ for $x\in N$ is an oriented
tangent $k$-plane. We say that $N$ is a {\it calibrated 
submanifold} if $\vp\vert_{T_xN}=\vol_{T_xN}$ for all~$x\in N$.
\label{uf2def1}
\end{dfn}

It is easy to show that calibrated submanifolds are automatically
{\it minimal submanifolds} \cite[Th.~II.4.2]{HaLa}. Here is the 
definition of SL $m$-folds in $\C^m$, taken from~\cite[\S III]{HaLa}.

\begin{dfn} Let $\C^m$ have complex coordinates $(z_1,\dots,z_m)$, 
and define a metric $g$, a real 2-form $\om$ and a complex $m$-form 
$\Om$ on $\C^m$ by
\begin{equation}
\begin{split}
g=\ms{\d z_1}+\cdots+\ms{\d z_m},\quad
\om&=\frac{i}{2}(\d z_1\w\d\bar z_1+\cdots+\d z_m\w\d\bar z_m),\\
\text{and}\quad\Om&=\d z_1\w\cdots\w\d z_m.
\end{split}
\label{uf2eq1}
\end{equation}
Then $\Re\Om$ and $\Im\Om$ are real $m$-forms on $\C^m$. Let
$L$ be an oriented real submanifold of $\C^m$ of real dimension 
$m$. We say that $L$ is a {\it special Lagrangian submanifold} 
of $\C^m,$ or {\it SL\/ $m$-fold}\/ for short, if $L$ is calibrated 
with respect to $\Re\Om$, in the sense of Definition~\ref{uf2def1}. 
\label{uf2def2}
\end{dfn}

Harvey and Lawson \cite[Cor.~III.1.11]{HaLa} give the following
alternative characterization of special Lagrangian submanifolds.

\begin{prop} Let\/ $L$ be a real\/ $m$-dimensional submanifold 
of\/ $\C^m$. Then $L$ admits an orientation making it into a
special Lagrangian submanifold of\/ $\C^m$ if and only if\/
$\om\vert_L\equiv 0$ and\/~$\Im\Om\vert_L\equiv 0$.
\label{uf2prop1}
\end{prop}

An $m$-dimensional submanifold $L$ in $\C^m$ is called {\it Lagrangian} 
if $\om\vert_L\equiv 0$. Thus special Lagrangian submanifolds are 
Lagrangian submanifolds satisfying the extra condition that 
$\Im\Om\vert_L\equiv 0$, which is how they get their name.

\subsection{Almost Calabi--Yau $m$-folds and SL $m$-folds} 
\label{uf22}

Probably the best general context for special Lagrangian
geometry is {\it almost Calabi--Yau manifolds}.

\begin{dfn} Let $m\ge 2$. An {\it almost Calabi--Yau $m$-fold}
is a quadruple $(X,J,\om,\Om)$ such that $(X,J)$ is a compact
$m$-dimensional complex manifold, $\om$ the K\"ahler form of
a K\"ahler metric $g$ on $X$, and $\Om$ a non-vanishing
holomorphic $(m,0)$-form on~$X$.

We call $(X,J,\om,\Om)$ a {\it Calabi--Yau $m$-fold} if in
addition
\begin{equation*}
\om^m/m!=(-1)^{m(m-1)/2}(i/2)^m\Om\w\bar\Om.
\end{equation*}
Then for each $x\in X$ there exists an isomorphism $T_xX\cong\C^m$
that identifies $g_x,\om_x$ and $\Om_x$ with the flat versions
$g,\om,\Om$ on $\C^m$ in \eq{uf2eq1}. Furthermore, $g$ is Ricci-flat
and its holonomy group is a subgroup of~$\SU(m)$.
\label{uf2def3}
\end{dfn}

This is not the usual definition of a Calabi--Yau manifold, but is 
essentially equivalent to it. Motivated by Proposition \ref{uf2prop1},
we define {\it special Lagrangian submanifolds} of almost Calabi--Yau
manifolds.

\begin{dfn} Let $(X,J,\om,\Om)$ be an almost Calabi--Yau $m$-fold 
with metric $g$, and $N$ a real $m$-dimensional submanifold of $X$. We 
call $N$ a {\it special Lagrangian submanifold}, or {\it SL $m$-fold} 
for short, if~$\om\vert_N\equiv\Im\Om\vert_N\equiv 0$. 
\label{uf2def4}
\end{dfn}

The properties of SL $m$-folds in almost Calabi--Yau $m$-folds
are discussed by the author in \cite{Joyc2}. It turns out
\cite[\S 9.5]{Joyc2} that SL $m$-folds in almost Calabi--Yau
$m$-folds are also calibrated w.r.t.\ $\Re\Om$, but using a
conformally rescaled metric~$\ti g=f^2g$.

The deformation and obstruction theory for {\it compact}\/ SL
$m$-folds in almost Calabi--Yau $m$-folds is well understood,
and beautifully behaved. Locally, SL $m$-folds in almost
Calabi--Yau $m$-folds are expected to behave like SL $m$-folds
in $\C^m$, especially in their singular behaviour. Thus, by
studying singular SL $m$-folds in $\C^m$, we learn about
singular SL $m$-folds in almost Calabi--Yau $m$-folds.

\section{Background material from analysis}
\label{uf3}

Here are some definitions we will need to make sense of
analytic results from \cite{Joyc3,Joyc4,Joyc5}. A closed,
bounded, contractible subset $S$ in $\R^n$ will be called
a {\it domain} if the {\it interior} $S^\circ$ of $S$ is
connected with $S=\overline{S^\circ}$, and the {\it boundary}
$\pd S=S\sm S^\circ$ is a compact embedded hypersurface in
$\R^n$. A domain $S$ in $\R^2$ is called {\it strictly convex}
if $S$ is convex and the curvature of $\pd S$ is nonzero at
every point.

Let $S$ be a domain in $\R^n$. Define $C^k(S)$ for $k\ge 0$ to be
the space of continuous functions $f:S\ra\R$ with $k$ continuous
derivatives, and $C^\iy(S)=\bigcap_{k=0}^\iy C^k(S)$. For $k\ge 0$
and $\al\in(0,1)$, define the {\it H\"older space} $C^{k,\al}(S)$
to be the subset of $f\in C^k(S)$ for which
\begin{equation*}
[\pd^k f]_\al=\sup_{x\ne y\in S}
\frac{\bmd{\pd^kf(x)-\pd^kf(y)}}{\md{x-y}^\al}
\quad\text{is finite.}
\end{equation*}

A {\it second-order quasilinear operator} $Q:C^2(S)\ra C^0(S)$
is an operator of the form
\begin{equation*}
\bigl(Qu\bigr)(x)=
\sum_{i,j=1}^na^{ij}(x,u,\pd u)\frac{\pd^2u}{\pd x_i\pd x_j}(x)
+b(x,u,\pd u),
\end{equation*}
where $a^{ij}$ and $b$ are continuous maps $S\t\R\t(\R^n)^*\ra\R$,
and $a^{ij}=a^{ji}$ for all $i,j=1,\ldots,n$. We call the functions
$a^{ij}$ and $b$ the {\it coefficients} of $Q$. We call $Q$ 
{\it elliptic} if the symmetric $n\t n$ matrix $(a^{ij})$ is 
positive definite at every point.

A second-order quasilinear operator $Q$ is in {\it divergence form}
if it is written
\begin{equation*}
\bigl(Qu\bigr)(x)=
\sum_{j=1}^n\frac{\pd}{\pd x_j}\bigl(a^j(x,u,\pd u)\bigr)
+b(x,u,\pd u)
\end{equation*}
for functions $a^j\in C^1\bigr(S\t\R\t(\R^n)^*\bigr)$ for
$j=1,\ldots,n$ and $b\in C^0\bigr(S\t\R\t(\R^n)^*\bigr)$.
If $Q$ is in divergence form, we say that integrable functions
$u,f$ are a {\it weak solution} of the equation $Qu=f$ if $u$
is weakly differentiable with weak derivative $\pd u$, and
$a^j(x,u,\pd u),b(x,u,\pd u)$ are integrable with
\begin{equation*}
-\sum_{j=1}^n\int_S\frac{\pd\psi}{\pd x_j}\cdot a^j(x,u,\pd u)
\d{\bf x}+\int_S\psi\cdot b(x,u,\pd u)\d{\bf x}
=\int_S\psi\cdot f\,\d{\bf x}
\end{equation*}
for all $\psi\in C^1(S)$ with~$\psi\vert_{\pd S}\equiv 0$.

If $Q$ is a second-order quasilinear operator, we may interpret
the equation $Qu=f$ in three different senses:
\begin{itemize}
\item We just say that $Qu=f$ if $u\in C^2(S)$, $f\in C^0(S)$
and $Qu=f$ in $C^0(S)$ in the usual way.
\item We say that $Qu=f$ {\it holds with weak derivatives} if $u$
is twice weakly differentiable and $Qu=f$ holds almost everywhere,
defining $Qu$ using weak derivatives.
\item We say that $Qu=f$ {\it holds weakly} if $Q$ is in divergence
form and $u$ is a weak solution of $Qu=f$. Note that this requires
only that $u$ be {\it once} weakly differentiable, and the second
derivatives of $u$ need not exist even weakly.
\end{itemize}

Clearly the first sense implies the second, which implies the third.
If $Q$ is {\it elliptic} and $a^j,b,f$ are suitably regular, one can
usually show that a weak solution to $Qu=f$ is a classical solution,
so that the three senses are equivalent. But for singular equations
that are not elliptic at every point, the three senses are distinct.

\section{Finding the equations}
\label{uf4}

Let $N$ be a special Lagrangian 3-fold in $\C^3$ invariant under
the $\U(1)$-action
\begin{equation}
{\rm e}^{i\th}:(z_1,z_2,z_3)\mapsto
({\rm e}^{i\th}z_1,{\rm e}^{-i\th}z_2,z_3)
\quad\text{for ${\rm e}^{i\th}\in\U(1)$.}
\label{uf4eq1}
\end{equation}
Locally we can write $N$ in the form
\begin{equation}
\begin{split}
N=\bigl\{(z_1,z_2,z_3)\in\C^3:\,& z_1z_2=v(x,y)+iy,\quad z_3=x+iu(x,y),\\
&\ms{z_1}-\ms{z_2}=2a,\quad (x,y)\in S\bigr\},
\end{split}
\label{uf4eq2}
\end{equation}
where $S$ is a domain in $\R^2$, $a\in\R$ and $u,v:S\ra\R$ are
continuous.

Here $\ms{z_1}-\ms{z_2}$ is the {\it moment map} of the $\U(1)$-action
\eq{uf4eq1}, and so $\ms{z_1}-\ms{z_2}$ is constant on any
$\U(1)$-invariant Lagrangian 3-fold in $\C^3$. We choose the constant
to be $2a$. Effectively \eq{uf4eq2} just means that we choose
$x=\Re(z_3)$ and $y=\Im(z_1z_2)$ as local coordinates on the
2-manifold $N/\U(1)$. Then we find~\cite[Prop.~4.1]{Joyc3}:

\begin{prop} Let\/ $S,a,u,v$ and\/ $N$ be as above. Then
\begin{itemize}
\item[{\rm(a)}] If\/ $a=0$, then $N$ is a (possibly singular)
special Lagrangian $3$-fold in $\C^3$ if\/ $u,v$ are
differentiable and satisfy
\begin{equation}
\frac{\pd u}{\pd x}=\frac{\pd v}{\pd y}
\quad\text{and}\quad
\frac{\pd v}{\pd x}=-2\bigl(v^2+y^2\bigr)^{1/2}\frac{\pd u}{\pd y},
\label{uf4eq3}
\end{equation}
except at points $(x,0)$ in $S$ with\/ $v(x,0)=0$, where $u,v$ 
need not be differentiable. The singular points of\/ $N$ are those
of the form $(0,0,z_3)$, where $z_3=x+iu(x,0)$ for $(x,0)\in S$ 
with\/~$v(x,0)=0$.
\item[{\rm(b)}] If\/ $a\ne 0$, then $N$ is a nonsingular special 
Lagrangian $3$-fold in $\C^3$ if and only if\/ $u,v$ are differentiable
in $S$ and satisfy
\begin{equation}
\frac{\pd u}{\pd x}=\frac{\pd v}{\pd y}\quad\text{and}\quad
\frac{\pd v}{\pd x}=-2\bigl(v^2+y^2+a^2\bigr)^{1/2}\frac{\pd u}{\pd y}.
\label{uf4eq4}
\end{equation}
\end{itemize}
\label{uf4prop1}
\end{prop}

The proof is elementary: at each point ${\bf z}\in N$ we calculate
the tangent space $T_{\bf z}N$ in terms of $\pd u,\pd v$, and use
Proposition \ref{uf2prop1} to find the conditions for $T_{\bf z}N$
to be a special Lagrangian $\R^3$ in $\C^3$. If ${\bf z}=(0,0,z_3)$
then $\d\bigl(\ms{z_1}-\ms{z_2}\bigr)=0$ at $\bf z$, so $\bf z$ is
a singular point of $N$, and $T_{\bf z}N$ does not exist.

Using \eq{uf4eq4} to write $\frac{\pd}{\pd y}\bigl(\frac{\pd u}{\pd x}
\bigr)$ and $\frac{\pd}{\pd x}\bigl(\frac{\pd u}{\pd y}\bigr)$ in
terms of $v$ and setting $\frac{\pd^{\smash{2}}u}{\pd y\pd x}=
\frac{\pd^{\smash{2}}u}{\pd x\pd y}$, we easily
prove~\cite[Prop.~8.1]{Joyc3}:

\begin{prop} Let\/ $S$ be a domain in $\R^2$ and\/ $u,v\in C^2(S)$
satisfy \eq{uf4eq4} for $a\ne 0$. Then
\begin{equation}
\frac{\pd}{\pd x}\Bigl[\bigl(v^2+y^2+a^2\bigr)^{-1/2}
\frac{\pd v}{\pd x}\Bigr]+2\,\frac{\pd^2v}{\pd y^2}=0.
\label{uf4eq5}
\end{equation}
Conversely, if\/ $v\in C^2(S)$ satisfies \eq{uf4eq5} then 
there exists $u\in C^2(S)$, unique up to addition of a 
constant\/ $u\mapsto u+c$, such that $u,v$ satisfy~\eq{uf4eq4}.
\label{uf4prop2}
\end{prop}

Now \eq{uf4eq5} is a second order quasilinear elliptic
equation, in divergence form. Thus we can consider
{\it weak solutions} of \eq{uf4eq5} when $a=0$, which need
be only once weakly differentiable. We shall be interested
in solutions of \eq{uf4eq3} with singularities, and the
corresponding SL 3-folds $N$. It will be helpful to define
a class of {\it singular solutions} of~\eq{uf4eq3}.

\begin{dfn} Let $S$ be a domain in $\R^2$ and $u,v\in C^0(S)$. We
say that $(u,v)$ is a {\it singular solution} of \eq{uf4eq3} if
\begin{itemize}
\item[(i)] $u,v$ are weakly differentiable, and their weak
derivatives $\frac{\pd u}{\pd x},\frac{\pd u}{\pd y},
\frac{\pd v}{\pd x},\frac{\pd v}{\pd y}$ satisfy~\eq{uf4eq3}.
\item[(ii)] $v$ is a {\it weak solution} of \eq{uf4eq5} with
$a=0$, as in~\S\ref{uf3}.
\item[(iii)] Define the {\it singular points} of $u,v$ to be
the $(x,0)\in S$ with $v(x,0)=0$. Then except at singular points,
$u,v$ are $C^2$ in $S$ and real analytic in $S^\circ$, and
satisfy \eq{uf4eq3} in the classical sense.
\item[(iv)] For $a\in(0,1]$ there exist $u_a,v_a\in C^2(S)$
satisfying \eq{uf4eq4} such that $u_a\ra u$ and $v_a\ra v$
in $C^0(S)$ as~$a\ra 0_+$.
\end{itemize}
\label{uf4def}
\end{dfn}

This list of properties is somewhat arbitrary. The point is that
\cite[\S 8--\S 9]{Joyc4} gives powerful existence and uniqueness
results for solutions $u,v$ of \eq{uf4eq3} satisfying conditions
(i)--(iv) and various boundary conditions on $\pd S$, and all of
(i)--(iv) are useful in different contexts.

\section{Examples}
\label{uf5}

Here are four examples of SL 3-folds $N$ in the form \eq{uf4eq2},
taken from~\cite[\S 5]{Joyc3}.

\begin{ex}
Let $a\ge 0$, and define
\begin{align*}
N_a=\Bigl\{(z_1&,z_2,z_3)\in\C^3:\ms{z_1}-2a=\ms{z_2}=\ms{z_3},\\
&\Im\bigl(z_1z_2z_3\bigr)=0,\quad \Re\bigl(z_1z_2z_3\bigr)\ge 0\Bigr\}.
\end{align*}
Then $N_a$ is a nonsingular SL 3-fold diffeomorphic to
${\mathcal S}^1\t\R^2$ when $a>0$, and $N_0$ is an SL $T^2$-cone
with one singular point at $(0,0,0)$. The $N_a$ are invariant
under the $\U(1)^2$-action
\begin{equation*}
({\rm e}^{i\th_1},{\rm e}^{i\th_2}):(z_1,z_2,z_3)\longmapsto
({\rm e}^{i\th_1}z_1,{\rm e}^{i\th_2}z_2,{\rm e}^{-i\th_1-i\th_2}z_3),
\end{equation*}
which includes the $\U(1)$-action \eq{uf4eq1}, and are part of a
family of explicit $\U(1)^2$-invariant SL 3-folds written down by
Harvey and Lawson \cite[\S III.3.A]{HaLa}. By \cite[Th.~5.1]{Joyc3},
these SL 3-folds can be written in the form \eq{uf4eq2}, using
functions~$u_a,v_a:\R^2\ra\R$.
\label{uf5ex1}
\end{ex}

\begin{ex} Let $\al,\be,\ga\in\R$ and define $u(x,y)=\al x+\be$ 
and $v(x,y)=\al y+\ga$. Then $u,v$ satisfy \eq{uf4eq4} for any 
value of~$a$. 
\label{uf5ex2}
\end{ex}

\begin{ex} Define $u(x,y)=y\tanh x$ and $v(x,y)=\ha y^2\sech^2x-
\ha\cosh^2x$. Then $u$ and $v$ satisfy \eq{uf4eq3}. Equation
\eq{uf4eq2} with $a=0$ defines an explicit nonsingular SL 3-fold
$N$ in $\C^3$. It arises from Harvey and Lawson's `austere
submanifold' construction \cite[\S III.3.C]{HaLa} of SL $m$-folds
in $\C^m$, as the normal bundle of a catenoid in~$\R^3$.
\label{uf5ex3}
\end{ex}

\begin{ex} Define $u(x,y)=\md{y}-\ha\cosh 2x$ and $v(x,y)=-y\sinh 2x$.
Then $u,v$ satisfy \eq{uf4eq3}, except that $\frac{\pd u}{\pd y}$ is 
not well-defined on the $x$-axis. So equation \eq{uf4eq2} with $a=0$
gives an explicit SL 3-fold $N$ in $\C^3$. It is the union of two
nonsingular SL 3-folds intersecting in a real curve, which are
constructed in \cite[Ex.~7.4]{Joyc1} by evolving paraboloids in~$\C^3$.
\label{uf5ex4}
\end{ex}

One can show that when $a=0$, all four examples yield
{\it singular solutions} of \eq{uf4eq3} in $\R^2$, in the
sense of Definition~\ref{uf4def}.

\section{Generating $u,v$ from a potential $f$}
\label{uf6}

If $u,v$ satisfy \eq{uf4eq4} then as $\frac{\pd u}{\pd x}=
\frac{\pd v}{\pd y}$ there exists a {\it potential} $f$ for
$u,v$ with $\frac{\pd f}{\pd y}=u$, $\frac{\pd f}{\pd x}=v$.
So we easily prove~\cite[Prop.~7.1]{Joyc3}:

\begin{prop} Let\/ $S$ be a domain in $\R^2$ and\/ $u,v\in C^1(S)$
satisfy \eq{uf4eq4} for $a\ne 0$. Then there exists $f\in C^2(S)$
with\/ $\frac{\pd f}{\pd y}=u$, $\frac{\pd f}{\pd x}=v$ and
\begin{equation}
\Bigl(\Bigl(\frac{\pd f}{\pd x}\Bigr)^2+y^2+a^2\Bigr)^{-1/2}
\frac{\pd^2f}{\pd x^2}+2\,\frac{\pd^2f}{\pd y^2}=0.
\label{uf6eq1}
\end{equation}
This $f$ is unique up to addition of a constant, $f\mapsto f+c$.
Conversely, all solutions of\/ \eq{uf6eq1} yield solutions 
of\/~\eq{uf4eq4}. 
\label{uf6prop1}
\end{prop}

Equation \eq{uf6eq1} is a second-order quasilinear elliptic equation,
singular when $a=0$, which may be written in divergence form. The
following condensation of \cite[Th.~7.6]{Joyc3} and \cite[Th.s 9.20
\& 9.21]{Joyc4} proves existence and uniqueness for the {\it Dirichlet
problem} for~\eq{uf6eq1}.

\begin{thm} Suppose $S$ is a strictly convex domain in $\R^2$ invariant
under $(x,y)\mapsto(x,-y)$, and\/ $k\ge 0$, $\al\in(0,1)$. Let\/
$a\in\R$ and\/ $\phi\in C^{k+3,\al}(\pd S)$. Then if\/ $a\ne 0$ there
exists a unique $f\in C^{k+3,\al}(S)$ with\/ $f\vert_{\pd S}=\phi$
satisfying \eq{uf6eq1}. If\/ $a=0$ there exists a unique $f\in C^1(S)$
with\/ $f\vert_{\pd S}=\phi$, which is twice weakly differentiable and
satisfies \eq{uf6eq1} with weak derivatives.

Define $u=\frac{\pd f}{\pd y}$ and\/ $v=\frac{\pd f}{\pd x}$. If\/
$a\ne 0$ then $u,v\in C^{k+2,\al}(S)$ satisfy \eq{uf4eq4}, and if\/
$a=0$ then $u,v\in C^0(S)$ are a singular solution of\/ \eq{uf4eq3},
in the sense of Definition $\ref{uf4def}$. Furthermore, $f$ depends
continuously in $C^1(S)$, and\/ $u,v$ depend continuously in $C^0(S)$,
on $(\phi,a)$ in~$C^{k+3,\al}(\pd S)\t\R$.
\label{uf6thm1}
\end{thm}

Here is a very brief sketch of the proof. As $f,u,v$ satisfy
certain linear elliptic equations, using the maximum principle
the maxima and minima of $f,u,v$ are achieved on $\pd S$. Thus
$\md{f}\le\sup_{\pd S}\md{\phi}$. Also, using linear functions
$f'=\al+\be x+\ga y$ as comparison solutions of \eq{uf6eq1} and
the strict convexity of $S$ we can bound $u,v$ on $\pd S$, and
hence on $S$, in terms of the first two derivatives of $\phi$.
So we have {\it a priori} bounds for $f$ in $C^1(S)$ and $u,v$
in~$C^0(S)$.

When $a\ne 0$, bounds on $v$ and $y$ imply that \eq{uf6eq1} is
{\it uniformly elliptic}. We can therefore use standard results
on the Dirichlet problem for second order, quasilinear, uniformly
elliptic equations to prove the existence of a unique $f\in
C^{k+3,\al}(S)$ satisfying \eq{uf6eq1} with~$f\vert_{\pd S}=\phi$.

However, when $a=0$ equation \eq{uf6eq1} is not uniformly elliptic,
and there do not appear to be standard results available to complete
the theorem. Therefore in \cite[\S 9]{Joyc4} we define $f_a\in
C^{k+3,\al}(S)$ for $a\in(0,1]$ to be the unique solution of
\eq{uf6eq1} with $f\vert_{\pd S}=\phi$, and we show that $f_a$
converges in $C^1(S)$ as $a\ra 0_+$ to a solution $f$ of
\eq{uf6eq1} for $a=0$, with weak second derivatives.

The key step in doing this is to prove {\it a priori} estimates
for $f_a$ and its first two derivatives, that hold uniformly for
$a\in(0,1]$. Proving such estimates, and making them strong enough
to ensure that $u,v$ are continuous, was responsible for the length
and technical difficulty of~\cite{Joyc4}.

Combining Proposition \ref{uf4prop1} and Theorem \ref{uf6thm1} gives
existence and uniqueness for a large class of $\U(1)$-invariant SL
3-folds in $\C^3$, with boundary conditions, including {\it singular}
SL 3-folds. It is interesting that this existence and uniqueness is
{\it entirely unaffected} by singularities appearing in~$S^\circ$. 

\section{Results motivated by complex analysis}
\label{uf7}

In \cite[\S 6]{Joyc3} and \cite[\S 7]{Joyc5} we study the zeroes
of $(u_1,v_1)-(u_2,v_2)$, where $(u_j,v_j)$ satisfy \eq{uf4eq3}
or \eq{uf4eq4}. A key tool is the idea of {\it winding number}.

\begin{dfn} Let $C$ be a compact oriented 1-manifold,
and $\ga:C\ra\R^2\sm\{0\}$ a differentiable map. Then
the {\it winding number of\/ $\ga$ about\/ $0$ along} $C$ 
is $\frac{1}{2\pi}\int_C\ga^*(\d\th)$, where $\d\th$ is
the closed 1-form $x^{-1}\d y-y^{-1}\d x$ on~$\R^2\sm\{0\}$.
\label{uf7def1}
\end{dfn}

The motivation is the following theorem from elementary
complex analysis:

\begin{thm} Let\/ $S$ be a domain in $\C$, and suppose
$f:S\ra\C$ is a holomorphic function, with\/ $f\ne 0$ on
$\pd S$. Then the number of zeroes of\/ $f$ in $S^\circ$,
counted with multiplicity, is equal to the winding number
of\/ $f\vert_{\pd S}$ about\/ $0$ along~$\pd S$.
\label{uf7thm1}
\end{thm}

As \eq{uf4eq3} and \eq{uf4eq4} are nonlinear versions of the
Cauchy--Riemann equations for holomorphic functions, it is
natural to expect that similar results should hold for their
solutions. The first step is to define the {\it multiplicity} of
an isolated zero in $S^\circ$, following~\cite[Def.~7.1]{Joyc5}.

\begin{dfn} Let $S$ be a domain in $\R^2$, and $a\in\R$. Suppose
$u_j,v_j:S\ra\R$ for $j=1,2$ are solutions of \eq{uf4eq4} in
$C^1(S)$ if $a\ne 0$, and singular solutions of \eq{uf4eq3} in
$C^0(S)$ if $a=0$, in the sense of Definition~\ref{uf4def}.

We call a point $(b,c)\in S$ a {\it zero} of $(u_1,v_1)-(u_2,v_2)$
in $S$ if $(u_1,v_1)=(u_2,v_2)$ at $(b,c)$. A zero $(b,c)$ is called
{\it singular} if $a=c=0$ and $v_1(b,0)=v_2(b,0)=0$, so that $(b,c)$
is a {\it singular point} of $u_1,v_1$ and $u_2,v_2$. Otherwise we
say $(b,c)$ is a {\it nonsingular zero}. We call a zero $(b,c)$
{\it isolated} if for some $\ep>0$ there exist no other zeroes $(x,y)$
of $(u_1,v_1)-(u_2,v_2)$ in $S$ with~$0<(x-b)^2+(y-c)^2\le\ep^2$.

Let $(b,c)\in S^\circ$ be an isolated zero of $(u_1,v_1)-(u_2,v_2)$.
Define the {\it multiplicity} of $(b,c)$ to be the winding number of
$(u_1,v_1)-(u_2,v_2)$ about 0 along the positively oriented circle
$\ga_\ep(b,c)$ of radius $\ep$ about $(b,c)$, where $\ep>0$ is chosen
small enough that $\ga_\ep(b,c)$ lies in $S^\circ$ and $(b,c)$ is the
only zero of $(u_1,v_1)-(u_2,v_2)$ inside or on~$\ga_\ep(b,c)$.
\label{uf7def2}
\end{dfn}

From \cite[\S 6.1]{Joyc3} and \cite[Cor.~7.6]{Joyc5}, we have:

\begin{thm} In the situation above, the multiplicity of any
isolated zero $(b,c)$ of\/ $(u_1,v_1)-(u_2,v_2)$ in $S^\circ$
is a positive integer.
\label{uf7thm2}
\end{thm}

The proof is different depending on whether $(b,c)$ is a
singular or a nonsingular zero. If $(b,c)$ is nonzero one
can show~\cite[Prop.~6.5]{Joyc3}:

\begin{prop} In the situation above, suppose $(u_1,v_1)-(u_2,v_2)$
has an isolated, nonsingular zero at $(b,c)\in S^\circ$. Then there
exists $k\ge 1$ and\/ $C\in \C\sm\{0\}$ such that
\begin{equation}
\begin{split}
\la u_1(x,y)+iv_1(x,y)=\la u_2(x,y)&+iv_2(x,y)
+C\bigl(\la(x-b)+i(y-c)\bigr)^k\\
&+O\bigl(\md{x-b}^{k+1}+\md{y-c}^{k+1}\bigr),
\end{split}
\label{uf7eq1}
\end{equation}
where~$\la=\sqrt{2}\bigl(v_1(b,c)^2+c^2+a^2\bigr)^{1/4}$.
\label{uf7prop1}
\end{prop}

The point is that $\la(u_1-u_2)+i(v_1-v_2)$ is like a holomorphic
function of $\la x+iy$ near $\la b+ic$, so to leading order it
is a multiple of $(\la(x-b)+i(y-c))^k$ for some $k\ge 1$. Comparing
\eq{uf7eq1} with Definition \ref{uf7def2} we see that $k$ is the
{\it multiplicity} of $(b,c)$, which proves that multiplicities
of nonsingular zeroes are positive integers.

Proposition \ref{uf7prop1} also gives an alternative, more
familiar characterization of the multiplicity of a nonsingular
zero, \cite[Def.~6.3]{Joyc3}: an isolated, nonsingular zero
$(b,c)$ of $(u_1,v_1)-(u_2,v_2)$ has multiplicity $k\ge 1$ if
$\pd^ju_1(b,c)=\pd^ju_2(b,c)$ and $\pd^jv_1(b,c)=\pd^jv_2(b,c)$
for $j=0,\ldots,k-1$, but $\pd^ku_1(b,c)\ne\pd^ku_2(b,c)$
and~$\pd^kv_1(b,c)\ne\pd^kv_2(b,c)$.

Proving Theorem \ref{uf7thm2} for singular zeroes is more tricky.
Basically it follows from the nonsingular case by a limiting
argument involving part (iv) of Definition \ref{uf4def}, but
there are subtleties in showing that a singular zero cannot
have multiplicity zero.

If $(u_1,v_1)\not\equiv(u_2,v_2)$ then all zeroes of
$(u_1,v_1)-(u_2,v_2)$ in $S^\circ$ are isolated,
\cite[Cor.~6.6]{Joyc3}, \cite[Th.~7.8]{Joyc5}.

\begin{thm} Let\/ $S$ be a domain in $\R^2$, and $a\in\R$. If\/
$a\ne 0$ let\/ $u_j,v_j\in C^1(S)$ satisfy \eq{uf4eq4} for $j=1,2$,
and if\/ $a=0$ let\/ $u_j,v_j\in C^0(S)$ be singular solutions of\/
\eq{uf4eq3} for $j=1,2$. Then either $(u_1,v_1)\equiv(u_2,v_2)$, or
there are at most countably many zeroes of\/ $(u_1,v_1)-(u_2,v_2)$
in $S^\circ$, all isolated.
\label{uf7thm3}
\end{thm}

To prove this, suppose $(u_1,v_1)\not\equiv(u_2,v_2)$ and $(b,c)$
is a nonsingular zero of $(u_1,v_1)-(u_2,v_2)$ in $S^\circ$. As
$u_j,v_j$ are real analytic near $(b,c)$, if $\pd^ju_1(b,c)=
\pd^ju_2(b,c)$ and $\pd^jv_1(b,c)=\pd^jv_2(b,c)$ for all $j\ge 1$
then $(u_1,v_1)\equiv(u_2,v_2)$, a contradiction. Hence there
exists a smallest $k\ge 1$ such that $\pd^ku_1(b,c)\ne\pd^ku_2(b,c)$
or $\pd^kv_1(b,c)\ne\pd^kv_2(b,c)$. One can then show following
Proposition \ref{uf7prop1} that \eq{uf7eq1} holds for some
$C\in\C\sm\{0\}$, and therefore $(b,c)$ is an isolated zero.

Hence, if $(u_1,v_1)\not\equiv(u_2,v_2)$ then nonsingular zeroes
of $(u_1,v_1)-(u_2,v_2)$ in $S^\circ$ are isolated. To prove that
singular zeroes are isolated requires a careful study of singular
solutions $(u,v)$ of \eq{uf4eq3} with $v=0$ on an interval on the
$x$-axis, carried out in~\cite[\S 6]{Joyc5}.

Following Theorem \ref{uf7thm1}, we can now prove
\cite[Th.~6.7]{Joyc3}, \cite[Th.~7.7]{Joyc5}:

\begin{thm} Let\/ $S$ be a domain in $\R^2$, and $a\in\R$. If\/
$a\ne 0$ let\/ $u_j,v_j\in C^1(S)$ satisfy \eq{uf4eq4} for $j=1,2$,
and if\/ $a=0$ let\/ $u_j,v_j\in C^0(S)$ be singular solutions of\/
\eq{uf4eq3} for $j=1,2$. Suppose $(u_1,v_1)\neq(u_2,v_2)$ at every
point of\/ $\pd S$. Then $(u_1,v_1)-(u_2,v_2)$ has finitely many
zeroes in $S$, all isolated. Let there be $n$ zeroes, with
multiplicities $k_1,\ldots,k_n$. Then the winding number of\/
$(u_1,v_1)-(u_2,v_2)$ about\/ $0$ along $\pd S$ is~$\sum_{i=1}^nk_i$.
\label{uf7thm4}
\end{thm}

Suppose $u_j,v_j$ come from a potential $f_j$ as in \S\ref{uf6},
with $f_j\vert_{\pd S}=\phi_j$. One can show directly that if the
winding number of $(u_1,v_1)-(u_2,v_2)$ about 0 along $\pd S$
is $k$, and $\phi_1-\phi_2$ has $l$ local maxima and $l$ local
minima on $\pd S$, then $\md{k}\le l-1$. So we prove
\cite[Th.~7.11]{Joyc3}, \cite[Th.~7.10]{Joyc5}:

\begin{thm} Suppose $S$ is a strictly convex domain in $\R^2$ invariant
under $(x,y)\mapsto(x,-y)$, and\/ $a\in\R$, $k\ge 0$, $\al\in(0,1)$,
and\/ $\phi_1,\phi_2\in C^{k+3,\al}(\pd S)$. Let\/ $u_j,v_j\in C^0(S)$
be the (singular) solution of\/ \eq{uf4eq3} or \eq{uf4eq4} constructed
in Theorem $\ref{uf6thm1}$ from $\phi_j$, for~$j=1,2$.

Suppose $\phi_1-\phi_2$ has $l$ local maxima and\/ $l$
local minima on $\pd S$. Then $(u_1,v_1)-(u_2,v_2)$ has finitely
many zeroes in $S^\circ$, all isolated. Let there be $n$ zeroes
in $S^\circ$ with multiplicities $k_1,\ldots,k_n$.
Then~$\sum_{i=1}^nk_i\le l-1$.
\label{uf7thm5}
\end{thm}

In particular, if $l=1$ then $(u_1,v_1)\ne(u_2,v_2)$ in $S^\circ$.
This will be useful in~\S\ref{uf9}.

\section{A rough classification of singular points}
\label{uf8}

We can now use \S\ref{uf7} to study singular points of $u,v$,
following~\cite[\S 9]{Joyc5}.

\begin{dfn} Let $S$ be a domain in $\R^2$, and $u,v\in C^0(S)$ a
singular solution of \eq{uf4eq3}, as in Definition \ref{uf4def}. Suppose
for simplicity that $S$ is invariant under $(x,y)\mapsto(x,-y)$.
Define $u',v'\in C^0(S)$ by $u'(x,y)=u(x,-y)$ and $v'(x,y)=-v(x,-y)$.
Then $u',v'$ is also a singular solution of~\eq{uf4eq3}.

A {\it singular point}, or {\it singularity}, of $(u,v)$ is a point
$(b,0)\in S$ with $v(b,0)=0$. Observe that a singularity of $(u,v)$
is automatically a zero of $(u,v)-(u',v')$. Conversely, a zero of
$(u,v)-(u',v')$ on the $x$-axis is a singularity. A singularity
of $(u,v)$ is called {\it isolated} if it is an isolated zero of
$(u,v)-(u',v')$. Define the {\it multiplicity} of an isolated
singularity $(b,0)$ of $(u,v)$ in $S^\circ$ to be the multiplicity
of $(u,v)-(u',v')$ at $(b,0)$, in the sense of Definition \ref{uf7def2}.
By Theorem \ref{uf7thm2}, this multiplicity is a positive integer.
\label{uf8def1}
\end{dfn}

From Theorem \ref{uf7thm3} we deduce~\cite[Th.~9.2]{Joyc5}:

\begin{thm} Let\/ $S$ be a domain in $\R^2$ invariant under
$(x,y)\mapsto(x,-y)$, and\/ $u,v\in C^0(S)$ a singular solution
of\/ \eq{uf4eq3}. If\/ $u(x,-y)\equiv u(x,y)$ and\/ $v(x,-y)\equiv
-v(x,y)$ then $(u,v)$ is singular along the $x$-axis in $S$, and
the singularities are nonisolated. Otherwise there are at most
countably many singularities of\/ $(u,v)$ in $S^\circ$, all isolated.
\label{uf8thm1}
\end{thm}

We divide isolated singularities $(b,0)$ into four types,
depending on the behaviour of $v(x,0)$ near~$(b,0)$.

\begin{dfn} Let $S$ be a domain in $\R^2$, and $u,v\in C^0(S)$ a
singular solution of \eq{uf4eq3}, as in Definition \ref{uf4def}.
Suppose $(b,0)$ is an isolated singular point of $(u,v)$ in $S^\circ$.
Then there exists $\ep>0$ such that $\,\overline{\!B}_\ep(b,0)
\subset S^\circ$ and $(b,0)$ is the only singularity of $(u,v)$
in $\,\overline{\!B}_\ep(b,0)$. Thus, for $0<\md{x-b}\le\ep$ we
have $(x,0)\in S^\circ$ and $v(x,0)\ne 0$. So by continuity $v$
is either positive or negative on each of $[b-\ep,b)\t\{0\}$
and~$(b,b+\ep]\t\{0\}$.
\begin{itemize}
\item[(i)] if $v(x)<0$ for $x\in[b-\ep,b)$ and $v(x)>0$ for
$x\in(b,b+\ep]$ we say the singularity $(b,0)$ is of
{\it increasing type}.
\item[(ii)] if $v(x)>0$ for $x\in[b-\ep,b)$ and $v(x)<0$ for
$x\in(b,b+\ep]$ we say the singularity $(b,0)$ is of
{\it decreasing type}.
\item[(iii)] if $v(x)<0$ for $x\in[b-\ep,b)$ and $v(x)<0$ for
$x\in(b,b+\ep]$ we say the singularity $(b,0)$ is of
{\it maximum type}.
\item[(iv)] if $v(x)>0$ for $x\in[b-\ep,b)$ and $v(x)>0$ for
$x\in(b,b+\ep]$ we say the singularity $(b,0)$ is of
{\it minimum type}.
\end{itemize}
\label{uf8def2}
\end{dfn}

The type determines whether the multiplicity is even or
odd,~\cite[Prop.~9.4]{Joyc5}.

\begin{prop} Let\/ $u,v\in C^0(S)$ be a singular solution of\/
\eq{uf4eq3} on a domain $S$ in $\R^2$, and\/ $(b,0)$ be an
isolated singularity of\/ $(u,v)$ in $S^\circ$ with multiplicity
$k$. If\/ $(b,0)$ is of increasing or decreasing type then $k$ is
odd, and if\/ $(b,0)$ is of maximum or minimum type then $k$ is even.
\label{uf8prop1}
\end{prop}

Theorem \ref{uf7thm5} yields a criterion for finitely many
singularities,~\cite[Th.~9.7]{Joyc5}:

\begin{thm} Suppose $S$ is a strictly convex domain in $\R^2$
invariant under $(x,y)\mapsto(x,-y)$, and\/ $\phi\in C^{k+3,\al}
(\pd S)$ for $k\ge 0$ and\/ $\al\in(0,1)$. Let\/ $u,v$ be the
singular solution of\/ \eq{uf4eq3} in $C^0(S)$ constructed
from $\phi$ in Theorem~$\ref{uf6thm1}$.

Define $\phi'\in C^{k+3,\al}(\pd S)$ by $\phi'(x,y)=-\phi(x,-y)$.
Suppose $\phi-\phi'$ has $l$ local maxima and\/ $l$ local minima
on $\pd S$. Then $(u,v)$ has finitely many singularities in
$S^\circ$. Let there be $n$ singularities in $S^\circ$ with
multiplicities $k_1,\ldots,k_n$. Then~$\sum_{i=1}^nk_i\le l-1$.
\label{uf8thm3}
\end{thm}

By applying Theorem \ref{uf6thm1} with $S$ the unit disc in $\R^2$
and $\phi$ a linear combination of functions $\sin(j\th),\cos(j\th)$
on the unit circle $\pd S$, we show~\cite[Cor.~10.10]{Joyc5}:

\begin{thm} There exist examples of singular solutions $u,v$
of\/ \eq{uf4eq3} with isolated singularities of every possible
multiplicity $n\ge 1$, and with both possible types allowed by
Proposition~$\ref{uf8prop1}$.
\label{uf8thm4}
\end{thm}

Combining this with Proposition \ref{uf4prop1} gives examples of
SL 3-folds in $\C^3$ with singularities of an {\it infinite number}
of different geometrical/topological types. We also show in
\cite[\S 10.4]{Joyc5} that singular points with multiplicity
$n\ge 1$ occur in {\it real codimension} $n$ in the family of
all SL 3-folds invariant under the $\U(1)$-action \eq{uf1eq1},
in a well-defined sense.

\section{Special Lagrangian fibrations}
\label{uf9}

We will now use our results to construct large families of
{\it special Lagrangian fibrations} of open subsets of $\C^3$
invariant under the $\U(1)$-action \eq{uf1eq1}, including
singular fibres. These will be important when we discuss the
{\it SYZ Conjecture} in \S\ref{uf10}, which concerns fibrations
of Calabi--Yau 3-folds by SL 3-folds.

\begin{dfn} Let $S$ be a strictly convex domain in $\R^2$ invariant
under $(x,y)\mapsto(x,-y)$, let $U$ be an open set in $\R^3$, and
$\al\in(0,1)$. Suppose $\Phi:U\ra C^{3,\al}(\pd S)$ is a continuous
map such that if $(a,b,c)\ne(a,b',c')$ in $U$ then $\Phi(a,b,c)-
\Phi(a,b',c')$ has exactly one local maximum and one local minimum
in~$\pd S$.

Let ${\bs\al}=(a,b,c)\in U$, and let $f_{\bs\al}\in
C^{3,\al}(S)$ be the unique (weak) solution of \eq{uf6eq1} with
$f_{\bs\al}\vert_{\pd S}=\Phi({\bs\al})$, which exists by Theorem
\ref{uf6thm1}. Define $u_{\bs\al}=\frac{\pd f_{\bs\al}}{\pd y}$
and $v_{\bs\al}=\frac{\pd f_{\bs\al}}{\pd x}$. Then
$(u_{\bs\al},v_{\bs\al})$ is a solution of \eq{uf4eq4} if
$a\ne 0$, and a singular solution of \eq{uf4eq3} if $a=0$.
Also $u_{\bs\al},v_{\bs\al}$ depend continuously on ${\bs\al}\in U$
in $C^0(S)$, by Theorem~\ref{uf6thm1}.

For each ${\bs\al}=(a,b,c)$ in $U$, define $N_{\bs\al}$ in $\C^3$ by
\begin{align*}
N_{\bs\al}=\bigl\{(z_1,z_2,z_3)\in\C^3:\,&
z_1z_2=v_{\bs\al}(x,y)+iy,\quad z_3=x+iu_{\bs\al}(x,y),\\
&\ms{z_1}-\ms{z_2}=2a,\quad (x,y)\in S^\circ\bigr\}.
\end{align*}
Then $N_{\bs\al}$ is a noncompact SL 3-fold without boundary in $\C^3$,
which is nonsingular if $a\ne 0$, by Proposition~\ref{uf4prop1}.
\label{uf9def1}
\end{dfn}

By \cite[Th.~8.2]{Joyc5} the $N_{\bs\al}$ are the fibres of an
{\it SL fibration}.

\begin{thm} In the situation of Definition $\ref{uf9def1}$, if\/
${\bs\al}\ne{\bs\al}'$ in $U$ then $N_{\bs\al}\cap N_{{\bs\al}'}
=\emptyset$. There exists an open set\/ $V\subset\C^3$ and a
continuous, surjective map $F:V\ra U$ such that\/ $F^{-1}({\bs\al})
=N_{\bs\al}$ for all\/ ${\bs\al}\in U$. Thus, $F$ is a special
Lagrangian fibration of\/ $V\subset\C^3$, which may include
singular fibres.
\label{uf9thm1}
\end{thm}

The main step in the proof is to show that distinct $N_{\bs\al}$
do not intersect, so that they fibre $V=\bigcup_{{\bs\al}\in U}
N_{\bs\al}$. Suppose $\al=(a,b,c)$ and $\al'=(a',b',c')$ are
distinct elements of $U$. If $a\ne a'$ then $N_{\bs\al}\cap
N_{{\bs\al}'}=\emptyset$, since $\ms{z_1}-\ms{z_2}$ is $2a$
on $N_{\bs\al}$ and $2a'$ on~$N_{{\bs\al}'}$.

If $a=a'$ then $\Phi({\bs\al})-\Phi({\bs\al}')$ has one local
maximum and one local minimum in $\pd S$, by Definition
\ref{uf9def1}. So Theorem \ref{uf7thm5} applies with $l=1$ to
show that $(u_{\bs\al},v_{\bs\al})-(u_{{\bs\al}'},v_{{\bs\al}'})$
has no zeroes in $S^\circ$, and again $N_{\bs\al}\cap N_{{\bs\al}'}
=\emptyset$. Thus distinct $N_{\bs\al}$ do not intersect.

Here is a simple way \cite[Ex.~8.3]{Joyc5} to produce families
$\Phi$ satisfying Definition \ref{uf9def1}, and thus generate
many SL fibrations of open subsets of~$\C^3$.

\begin{ex} Let $S$ be a strictly convex domain in $\R^2$ invariant
under $(x,y)\mapsto(x,-y)$, let $\al\in(0,1)$ and $\phi\in C^{3,\al}
(\pd S)$. Define $U=\R^3$ and $\Phi:\R^3\ra C^{3,\al}(\pd S)$ by
$\Phi(a,b,c)=\phi+bx+cy$. If $(a,b,c)\ne(a,b',c')$ then
$\Phi(a,b,c)-\Phi(a,b',c')=(b-b')x+(c-c')y\in C^\iy(\pd S)$. As
$b-b',c-c'$ are not both zero and $S$ is strictly convex, it
easily follows that $(b-b')x+(c-c')y$ has one local maximum
and one local minimum in~$\pd S$.

Hence the conditions of Definition \ref{uf9def1} hold for $S,U$
and $\Phi$, and so Theorem \ref{uf9thm1} defines an open set
$V\subset\C^3$ and a special Lagrangian fibration $F:V\ra\C^3$.
One can also show that changing the parameter $c$ in $U=\R^3$
just translates the fibres $N_{\bs\al}$ in $\C^3$, and
\begin{equation*}
V=\bigl\{(z_1,z_2,z_3)\in\C^3:(\Re z_3,\Im z_1z_2)\in S^\circ\bigr\}.
\end{equation*}
\label{uf9ex1}
\end{ex}

Here is very explicit example, taken from~\cite{Joyc6}.

\begin{ex} Define $F:\C^3\ra\R\t\C$ by
\begin{gather*}
F(z_1,z_2,z_3)=(a,b),\quad\text{where}\quad 2a=\ms{z_1}-\ms{z_2} \\
\text{and}\quad
b=\begin{cases}
z_3, & a=z_1=z_2=0, \\
z_3+\bar z_1\bar z_2/\md{z_1}, & a\ge 0,\;\> z_1\ne 0, \\
z_3+\bar z_1\bar z_2/\md{z_2}, & a<0.
\end{cases}
\end{gather*}
This is a piecewise-smooth SL fibration of $\C^3$. It is not
smooth on~$\md{z_1}=\md{z_2}$.

The fibres $F^{-1}(a,b)$ are $T^2$-cones singular at $(0,0,b)$ 
when $a=0$, and nonsingular ${\mathcal S}^1\t\R^2$ when $a\ne 0$. 
They are isomorphic to the SL 3-folds $N_{\md{a}}$ of Example
\ref{uf5ex1} under transformations of $\C^3$, but they are
assembled to make a fibration in a novel way. As $a$ goes
from positive to negative the fibres undergo a surgery, a
Dehn twist on~${\mathcal S}^1$.
\label{uf9ex2}
\end{ex}

\section{The SYZ Conjecture}
\label{uf10}

{\it Mirror Symmetry} is a mysterious relationship between pairs 
of Calabi--Yau 3-folds $X,\hat X$, arising from a branch of 
physics known as {\it String Theory}, and leading to some very 
strange and exciting conjectures about Calabi--Yau 3-folds, 
many of which have been proved in special cases.

Roughly speaking, String Theorists believe that each Calabi--Yau
3-fold $X$ has a quantization, a {\it Super Conformal Field
Theory} (SCFT). Invariants of $X$ such as the Dolbeault groups
$H^{p,q}(X)$ and the number of holomorphic curves in $X$
translate to properties of the SCFT. However, different
Calabi--Yau 3-folds $X,\hat X$ may have the same SCFT.

One way for this to happen is for the SCFT's of $X,\hat X$ to
be related by a certain simple involution of SCFT structure,
which does {\it not\/} correspond to a classical automorphism
of Calabi--Yau 3-folds. We then say that $X$ and $\hat X$ are
{\it mirror} Calabi--Yau 3-folds. One can argue using String
Theory that $H^{1,1}(X)\cong H^{2,1}(\hat X)$ and $H^{2,1}(X)
\cong H^{1,1}(\hat X)$. The mirror transform also exchanges
things related to the complex structure of $X$ with things
related to the symplectic structure of $\hat X$, and vice versa.

The {\it SYZ Conjecture}, due to Strominger, Yau and Zaslow
\cite{SYZ} in 1996, gives a geometric explanation of Mirror
Symmetry. Here is an attempt to state it.
\medskip

\noindent{\bf The SYZ Conjecture.} {\it Suppose $X$ and\/ $\hat X$ are 
mirror Calabi--Yau $3$-folds. Then (under some additional conditions) 
there should exist a compact topological\/ $3$-manifold\/ $B$ and 
surjective, continuous maps $f:X\ra B$ and\/ $\hat f:\hat X\ra B$, 
such that
\begin{itemize}
\item[{\rm(i)}] There exists a dense open set\/ $B_0\subset B$, such 
that for each\/ $b\in B_0$, the fibres $f^{-1}(b)$ and\/ $\hat f^{-1}(b)$
are nonsingular special Lagrangian $3$-tori $T^3$ in $X$ and\/ $\hat X$.
Furthermore, $f^{-1}(b)$ and\/ $\hat f^{-1}(b)$ are in some sense
dual to one another.
\item[{\rm(ii)}] For each\/ $b\in \De=B\sm B_0$, the fibres $f^{-1}(b)$ 
and\/ $\hat f^{-1}(b)$ are expected to be singular special Lagrangian
$3$-folds in $X$ and\/~$\hat X$.
\end{itemize}}
\medskip

We call $f,\hat f$ {\it special Lagrangian fibrations}, and the
set of singular fibres $\De$ is called the {\it discriminant}. It
is not yet clear what the final form of the SYZ Conjecture should
be: there are problems to do with the singular fibres, and with
what extra conditions on $X,\hat X$ are needed for $f,\hat f$
to exist.

Much mathematical research on the SYZ Conjecture has simplified
the problem by supposing that $f,\hat f$ are {\it Lagrangian
fibrations}, making only limited use of the `special' condition,
and supposing in addition that $f,\hat f$ are {\it smooth} maps.
Gross \cite{Gros1,Gros2,Gros3,Gros4}, Ruan \cite{Ruan1,Ruan2},
and others have built up a beautiful, detailed picture of how
dual SYZ fibrations work at the level of global symplectic
topology, in particular for examples such as the quintic and
its mirror, and for Calabi--Yau 3-folds constructed as
hypersurfaces in toric 4-folds, using combinatorial data.

The author's approach to the SYZ Conjecture \cite{Joyc6} has
a different viewpoint, and more modest aims. We take the
special Lagrangian condition seriously from the outset, and
focus on the local behaviour of SL fibrations near singular
points, rather than on global topological questions. Also,
we are interested in {\it generic} SL fibrations.

The best way to introduce a genericity condition is to
consider SL fibrations $f:X\ra B$ in which $X$ is a generic
almost Calabi--Yau 3-fold. The point of allowing $X$ to be
an {\it almost\/} Calabi--Yau 3-fold is that the family of
almost Calabi--Yau structures is infinite-dimensional, and
so picking a generic one is a strong condition, and should
simplify the singular behaviour of~$f$.

Now from \S\ref{uf9} we know a lot about $\U(1)$-invariant
SL fibrations of subsets of $\C^3$. By considering when these
are appropriate local models for singularities of SL fibrations
of almost Calabi--Yau 3-folds, the author makes the following
tentative suggestions:

\begin{itemize}
\item In a generic SL fibration $f:X\ra B$, the singularities
of codimension 1 in $B$ are locally modelled on the explicit
SL fibration $F$ given in Example~\ref{uf9ex2}.
\item Similarly, in generic SL fibrations $f:X\ra B$, one kind
of singular behaviour of codimension 2 in $B$ is modelled on a
$\U(1)$-invariant SL fibration of the kind considered in
\S\ref{uf9}, including a 1-parameter family of singular fibres
with isolated singular points of multiplicity 2, in the sense
of Definition~\ref{uf8def1}.
\item However, I do not expect codimension 3 singularities in
generic SL fibrations to be locally $\U(1)$-invariant, so this
approach will not help.
\end{itemize}

Here are some broader conclusions, also conjectural.

\begin{itemize}
\item For generic almost Calabi--Yau 3-folds $X$, SL fibrations
$f:X\ra B$ will not be smooth maps, but only {\it piecewise smooth}.
\item In a generic SL fibration $f:X\ra B$ the discriminant $\De$
is of codimension 1 in $B$, and the singular fibres are singular
at finitely many points.

In contrast, in the smooth Lagrangian fibrations $f:X\ra B$
considered by Gross and Ruan, the discriminant $\De$ is of
codimension 2 in $B$, and the typical singular fibre is singular
along an~${\mathcal S}^1$.
\item If $X,\hat X$ are a mirror pair of generic (almost)
Calabi--Yau 3-folds and $f:X\ra B$ and $\hat f:\hat X\ra B$ are
dual SL fibrations, then in general the discriminants $\De$ of $f$
and $\hat\De$ of $\hat f$ will {\it not\/} coincide in $B$. This
contradicts part (ii) of the SYZ Conjecture, as stated above.
\end{itemize}

For more details, see \cite{Joyc6}. In the author's view, these
calculations support the idea that the SYZ Conjecture in its
present form should be viewed primarily as a limiting statement,
about what happens at the `large complex structure limit', rather
than as simply being about pairs of Calabi--Yau 3-folds. A similar
conclusion is reached by Mark Gross in~\cite[\S 5]{Gros4}.

\end{document}